\documentclass{amsart}

\usepackage{graphicx}
\usepackage{multirow}
\usepackage{amsmath,amssymb,amsfonts}
\usepackage{mathrsfs}
\usepackage[title]{appendix}
\usepackage{xcolor}
\usepackage{textcomp}
\usepackage{manyfoot}
\usepackage{booktabs}
\usepackage{algorithm}
\usepackage{algorithmicx}
\usepackage{algpseudocode}
\usepackage{listings,hyperref}

\newtheorem{theorem}{Theorem}[section]

\newtheorem{lemma}[theorem]{Lemma}
\newtheorem{corollary}[theorem]{Corollary}

\newtheorem{example}[theorem]{Example}
\newtheorem{definition}[theorem]{Definition}
\newtheorem*{theorema}{Theorem}
\newcommand{\NN}{\mathbb N}
\newcommand{\ZZ}{\mathbb Z}
\newcommand{\RR}{\mathbb R}

\title[Gap estimates for the spectrum of the $m$-bonacci numbers]{Gap estimates for the spectrum of $m$-bonacci numbers}

\author{Anna Chiara Lai}
\address{
Dipartimento di Scienze di Base e Applicate per l'Ingegneria,
Sapienza Universit\`a di Roma,
Via Antonio Scarpa 16,
00161 Roma, Italy}
\email{annachiara.lai@uniroma1.it}

\author{Paola Loreti}
\address{
Dipartimento di Scienze di Base e Applicate per l'Ingegneria,
Sapienza Universit\`a di Roma,
Via Antonio Scarpa 16,
00161 Roma, Italy}
\email{paola.loreti@uniroma1.it}



\begin{document}

\begin{abstract}
We establish some results on the structure of the spectrum of $m$-bonacci numbers. More precisely, we study explicit lower bounds for the distance between elements separated by $N$ positions in the ordered spectrum. The result is further detailed in the Fibonacci and Tribonacci case. The methodology combines the combinatorial structure of $m$-bonacci words and the canonical $m$-bonacci number system.
\end{abstract}

\maketitle

\section{Introduction}
In this paper we are interested in gap estimates in the spectrum of $m$-bonacci numbers. 

To illustrate our results, we begin by considering the spectrum $\Lambda (\varphi)$ of the golden mean $\varphi$, namely the set of real numbers of the form $\sum_{i=0}^k \varepsilon_i \varphi^i$
with $k\in \mathbb{N}$ and $\varepsilon_i \in \{0,1\}$. Ordering $\Lambda(\varphi)$ strictly increasingly, so that $\Lambda(\varphi)=\{0<\lambda_1(\varphi)<\cdots<\lambda_n(\varphi)<\cdots\}$ we get
$$\Lambda(\varphi)=\{0,1,\varphi,\varphi^2,\varphi^2+1,\varphi^2+\varphi,\varphi^2+\varphi+1,\dots\}$$ and the consecutive distances (or gaps) $\lambda_{n+1}(\varphi)-\lambda_{n}(\varphi)$ take only two possible values $1$ and $\varphi-1$. If one considers larger increments $\lambda_{n+2}(\varphi)-\lambda_{n}(\varphi)$, the possible values are $2$ and $\varphi$. In particular, there is no $n\in\NN$ such that $\lambda_{n+2}(\varphi)-\lambda_{n}(\varphi)=2\varphi-2$. This constraint can be clearly understood by coding the gap sequence $\lambda_{n+1}(\varphi)-\lambda_n(\varphi)$ by assigning the symbol $1$ or $2$ depending on whether the gap is equal to $1$ or $\varphi-1$. The resulting sequence is $121121\cdots$, which coincides with the Fibonacci word. The Fibonacci word can be defined as the limit of the substitution $1 \mapsto 12$, $2 \mapsto 1$, obtained by iterating this rule starting from the symbol $1$ -- see \cite{AS03}.  In particular, the subword $22$ never occurs; this excludes the possibility of having two consecutive gaps equal to $\varphi-1$, and consequently $\lambda_{n+2}(\varphi)-\lambda_n(\varphi) \neq 2\varphi - 2$ for all $n$. More generally, finer combinatorial properties of the Fibonacci word,
such as balance and forbidden patterns, yield gap estimates of the form
\begin{equation}\label{lambdaest}\lambda_{n+N}(\varphi)-\lambda_{n}(\varphi)\geq N\left(\frac{1}{\varphi}+\frac{1}{\varphi^3}\right) -\left(\frac{1}{\varphi^4}+\varphi\right)\end{equation}
for all $N\geq 2,\,n\geq 1$ -- see Corollary \ref{cor1}. Writing the lower bound in \eqref{lambdaest} in the form $N\gamma_{2,N}$ with $\gamma_{2,N}:=\left(\frac{1}{\varphi}+\frac{1}{\varphi^3}\right) -\frac{1}{N}\left(\frac{1}{\varphi^4}+\varphi\right)$, we obtain that $\gamma_{2,N}\to \frac{1}{\varphi}+\frac{1}{\varphi^3}=\frac{\sqrt{5}}{\varphi}=\frac{1}{D^+(\Lambda(\varphi))}$ as $N\to \infty$, where $D^+(\Lambda(\varphi))$ denotes the (upper) density of $\Lambda(\varphi)$  -- see the definitions below. 

Our goal is to establish estimates of the form \eqref{lambdaest} for the wider class of $m$-bonacci numbers $q_m$, defined as the largest roots of the polynomials 
$$x^m-x^{m-1}-\cdots -x-1 \qquad \text{for } m\geq 2.$$
Note that the case $m=2$ corresponds to the golden mean case, because $q_2=\varphi$. We consider the spectrum of the $m$-bonacci numbers
$$\Lambda(q_m):=\left\{\sum_{i=0}^n\varepsilon_iq_m^i\mid \varepsilon_{i}\in\{0,1\},\, n=0,1,\dots\right\}=\{0<\lambda_1(q_m)<\lambda_2(q_m)<\cdots\}$$
and we establish gap estimates  of the form 
$$\lambda_{n+N}(q_m)-\lambda_{n}(q_m)\geq N\gamma_{m,N} \quad \text{for all }n\geq 1$$
 for any $N\geq 1$. Our main result provides a constant $\gamma_{m,N}$ which can be interpreted as a weighted average of the gap values in $\Lambda(q_m)$. The weights depend on the $m$-bonacci expansion of $N$,  and on the combinatorial structure of
the $m$-bonacci word coding the gaps. 
An $m$-bonacci expansion is a representation
of $N$ as a sum 
 of $m$-bonacci numbers $F^{(m)}_n$ with coefficients in $\{0,1\}$. The $m$-bonacci numbers are recursively defined by 
$F^{(m)}_0=1$ $F^{(m)}_1=2, \cdots F^{(m)}_{m-1}=2^{m-1}$
and, for $n\geq m$,
$$F^{(m)}_n=F^{(m)}_{n-1}+\cdots +F^{(m)}_{n-m}.$$
On the other hand, we recall that the $m$-bonacci word is defined as the fixed point of the Rauzy
substitution $\sigma_m$ defined by 
$$1\mapsto 12,\quad 2\mapsto 13,\quad \dots, \quad m-1\mapsto 1m, \quad m\mapsto 1.$$
In particular, the weights involve the balance constant $b_m$, which measures how uniformly the letters are distributed the $m$-bonacci word. 
Our main result is the following.  
\begin{theorem}\label{thm1i}
Let $m\geq 2$ and $N\in \mathbb{N}$. Then there exists a constant
$\gamma_{m,N}>0$ such that
\begin{equation}\label{gamma}
\lambda_{n+N}(q_m)-\lambda_n(q_m)\geq N\,\gamma_{m,N}
\quad \text{for all } n\geq 1.
\end{equation}

The following choice of $\gamma_{m,N}$ satisfies the above inequality.
Let $(x_h)\in\{0,1\}^\infty$ be the canonical $m$-bonacci expansion of $N$, that is,
a binary sequence with no occurrences of the block $1^m$ that expands $N$ in the form
\[
N=\sum_{h\geq 0} x_h F^{(m)}_h,
\quad x_h\in\{0,1\},
\]
where $(F^{(m)}_h)$ denotes the $m$-bonacci sequence.

Then
\begin{equation}\label{eqgamma}
\gamma_{m,N}
=\frac{1}{N}
\sum_{j=1}^m
\left(
\sum_{h\geq 0} x_h F^{(m)}_{h-j}
- b_m
\right)
d_m(j),
\end{equation}
where:
\begin{itemize}
\item[-] $d_m(j):=q_m^{j-1}-\sum_{k=0}^{j-2} q_m^k$ with $j=1,\dots,m$ are the finitely many possible gap values in $\Lambda(q_m)$;
\item[-] $b_m$ is a constant lower or equal than $m-1$ measuring the balance of the $m$-bonacci word.
\end{itemize}
\end{theorem}
We remark that explicit values for $b_m$ can be assured by results on the balance of $m$-bonacci word, borrowing from \cite{bal,balancetrib,generalbalance} we have that \eqref{eqgamma} holds by setting 
\begin{equation}\label{bm}
b_m:=\begin{cases}
    m-1 & \text{for } m=2,3,4;\\
    \lceil \frac{m+1}{2}\rceil& \text{for } m=5,\dots,12;\\
    \lfloor \kappa m\rfloor+12 & \text{for } m>12
\end{cases}\end{equation}
where $\kappa:= \frac{2}{\pi} \int_{0}^{2\pi}
\frac{1-\cos x}{(5-4\cos x)\ln(5-4\cos x)} \, dx \approx 0.58.$

By applying the above result to $m=2$ and remarking that the balance constant in the Fibonacci word is $b_2=1$, we obtain in the golden mean case 
$$\lambda_{N+n}(\varphi)-\lambda_{n}(\varphi)\geq \gamma_{2,N}=
\frac{1}{N}
\left(
\sum_{h=0}^\infty x_h(\varphi F_{h-1}+F_{h-2})
-\varphi^2
\right).$$
where $(x_h)\in\{0,1\}^\infty$ is the canonical Fibonacci expansion of $N$. Note that, in general, $m$-bonacci expansions $(x_h)$ 
have all $x_h=0$ except for finitely many. In Corollary \ref{cor1} and in Corollary \ref{cor2} we discuss in detail the cases $m=2,3$, respectively. 
\medskip

Let us place our contribution in the context of the literature on the spectra associated with Pisot numbers and gap distribution. 
The paper \cite{EJK90} is among the first to investigate the spectrum of a real number $q>1$, defined as:
$$\Lambda(q):=\left\{\sum_{i=0}^n\varepsilon_iq^i\mid \varepsilon_{i}\in\{0,1\},\, n=0,1,\dots\right\}=\{0<\lambda_1(q)<\lambda_2(q)<\cdots\}.$$
 Erd{\"o}s, Jo{\'o}, and Komornik were interested in the gap size $\lambda_{n+1}(q)-\lambda_n(q)$ showing that it is uniformly bounded from above by $1$ for all $q\in (1,2)$. The interest then moved to lower estimates: in the same paper it is proved that if $\lambda_{n+1}(q)-\lambda_n(q)\to 0$   then $1$ has an infinite expansion containing
arbitrarily long sequences of consecutive $0$ digits, namely $1=\sum_{i=1}^\infty q^{n_i}$ for some diverging sequence $n_i$ of natural numbers satisfying $\limsup_{i\to\infty} (n_{i+1}-n_i)=+\infty$. 
These results paved the way to many other investigations. In \cite{EJ92}, Erd{\"o}s, Jo{\'o}, and Jo{\'o}  showed that the $m$-bonacci number  $q_m$  is the unique root of $x^m-x^{m-1}-\cdots-x-1$ in $[1,2]$, $q_m\to2$ monotone increasingly as $m\to \infty$. They considered 
the finite sets
$\Lambda^{(h)}(q):=\left\{\sum_{k=1}^h \varepsilon_iq^i\mid \varepsilon_i\in\{0,1\}\right\}$
  for $h\geq 1$ and proved that the minimum gap in $\Lambda^{(h)}(q_m)$ is greater or equal to $\frac{1}{q_m}$ and equality holds for $h\geq m+1$.  As remarked in \cite{EJK98}, this result can be interpreted as follows. Setting for all $q>0$
  $$\ell(q):=\liminf_{n} \lambda_{n+1}(q)-\lambda_{n}(q)$$
  we obtain 
  $$\ell(q_m)=\frac{1}{q_m} \quad \text{for all }m\geq 2. $$
  In \cite{bug96}, Bugeaud considers the wider spectrum \begin{equation}\label{genspectrum}\Lambda^r(q):=\left\{\sum_{i=0}^n\varepsilon_iq^i\mid \varepsilon_{i}\in\{0,1,\dots,r\},\, n=0,1,\dots\right\}=\{0<\lambda_1^r(q)<\cdots<\lambda_n^r(q)\}\end{equation}
and the related quantity 
$$\ell^r(q):=\liminf_{n} (\lambda^r_{n+1}(q)-\lambda^r_n(q)).$$
Bugeaud shows that a real number $q\in(1,2)$ is a \emph{Pisot number}, namely a real algebraic integer  whose algebraic conjugates have modulus stricly less than $1$, if and only if $\ell^r(q)>0$ for all $r\geq 1$. Note that, for $m\geq 2$, $q_m$ is a Pisot number so Bugeaud's result applies. In \cite{EJK98} it is  proved that if $q>1$ is Pisot, then $\ell(q)>0$.
In the paper \cite{KLP00}, Komornik, Loreti, and Pedicini provide a complete description of $\ell^r(\varphi)$ when $\varphi$ is the golden mean: if $F_k$ is the $k$-th Fibonacci number $F_0=0,F_1=1,\, F_{k}=F_{k-1}+F_{k-2}$, and $\varphi^{k-2}< r\leq \varphi^{k-1}$ then $\ell^r(\varphi)=|F_{k}\varphi-F_{k+1}|$. 

Now, in \cite{FW02}, Feng and Wen investigate the structure of the set $\Lambda(q)$ in the Pisot case, showing that the gap sequence $\lambda_{n+1}(q)-\lambda_n(q)$ admits a finite number of values and that it can be coded via substitutions on a finite alphabet. A substitution on a finite alphabet assigns to each symbol a finite word over the same alphabet, and is then extended to longer words by concatenation. As mentioned above, by \cite{FW02,bug02}, 
when $q=q_m$ is the largest root of $x^m-x^{m-1}-\cdots-x-1$, the gaps in $\Lambda(q_m)$ are related to the Rauzy substitution $\sigma_m$
via its fixed point 
$$\mathbf v_{m}=\mathbf v_{m,1}\mathbf v_{m,2}\cdots:=\lim_{n\to\infty} \sigma_m^n(1)$$
that encodes the combinatorial structure of the spectrum. Explicit gap formulae are given in \cite{bug02} and, as corollary, Bugeaud proved in the same paper that the gap $1$ occurs in $\Lambda(q_m)$ with frequency $\frac{1}{q_m}$ whereas the gap 
$q_m^{v-1}-q_m^{v-2}-\cdots q_m-1$ occurs with frequency $\frac{1}{q_m^{v}}$. In \cite{H04} Hare gives an algorithm for the computation of the frequencies of gap sizes in more general cases -- this result is generalized to a wider classes of Pisot numbers in the paper \cite{GH06} by Garth and Hare. Akyiama and Komornik \cite{AK13} prove that the set generalized spectrum $Y^r(q):=\{\sum_{i=0}^n\varepsilon_i q^i\mid \varepsilon_i\in\{0,\pm 1,\dots,\pm r\}\}$ has not accumulation points if and only if $q$ is a Pisot number or $q\geq m+1$.   More recently, Feng  shows in \cite{F16} that for all $r\geq 1$ we have that $\ell^r(q)>0$ if and only if either $q$ is Pisot or $r<q-1$. 

Finally, we mention that the spectrum $\Lambda(q)$ of a Pisot number
is closely related to several fundamental objects in aperiodic order
and numeration theory. To illustrate this connection, we note that
$\Lambda(q_m)$ and the $m$-bonacci quasicrystal can be described in terms of the same substitution, although they give rise to different sets of gap sizes. For an overview, we refer to the survey by Berthé and Yassawi \cite{BY24} and the references therein.
\medskip

\emph{Organization of the paper.}
In Section \ref{sec2} we collect the preliminary material needed for the proof,
with particular emphasis on the combinatorial structure of the
$m$-bonacci word, and gap frequencies.
In Section \ref{sec3} we prove Theorem~\ref{thm1i}, combining the
$m$-bonacci expansions of integers with the balance properties of the $m$-bonacci word,  and we derive further estimates  for the Fibonacci case and the Tribonacci case. 
In Section \ref{sec4} we explore some density properties of
$\Lambda(q_m)$. 
Finally, in Section \ref{sec5} we present our conclusions.
\section{Preliminaries}\label{sec2}
We recall the definitions and results on $m$-bonacci words, $m$-bonacci expansions and density of $\Lambda(q_m)$.

\subsection{The $m$-bonacci words and their relation with $\Lambda(q_m)$}
For any $m\geq 2$, the $m$-bonacci word $\mathbf v_m$ is generated iteratively, via the repeated application of a substitution map $\sigma_m$ --known as the \emph{Rauzy substitution}-- up to convergence to its fixed point, which is $\mathbf v_m$. 

\begin{definition}[Rauzy substitution and $m$-bonacci word]   Let $\mathcal A_m:=\{1,\dots,m\}$ with $m\geq 2$ and $\mathcal A_m^*$ be the set of finite words with digits in $\mathcal A_m$. For every $m\geq 2$, we consider \emph{Rauzy substitution} $\sigma_m:\mathcal A_m^*\to \mathcal A_m^*$ given by
$k\mapsto 1(k+1)$ for $k=1,\dots,m-1$ and $m\mapsto 1$.

For all $m\geq 2$ let
    $$\mathbf v_m:=\lim_{n\to \infty}\sigma^n_m(1).$$ 
    
   The words $\mathbf v_2=\sigma_2^\infty(1)=121121\cdots$  and $\mathbf v_3=\sigma_3^\infty(1)=121312\cdots$ are respectively called the \emph{Fibonacci word} and \emph{Tribonacci word}. For $m> 3$, the word $\mathbf v_m$ is called the \emph{$m$-bonacci word}. 
\end{definition}
For a general introduction on substitutions and, in particular, on the well-posedness of the limit in the definition of $\mathbf v_m$ we refer to \cite{Loth}. Here we limit ourselves to remarking that the substitution $\sigma_m(1)$ begins with $1$. 
Hence the sequence of words $(\sigma_m^n(1))_{n\geq 1}$ is increasing
for the prefix order and thus converges to a unique infinite word.

\medskip
The next result relates $m$-bonacci words to the gaps in the spectrum $\Lambda(q_m)$ in terms of size and frequency. Recall that a gap $g$ occurs with frequency $f$ in an increasingly ordered set $\{\lambda_k\}$ if 
$$\lim_{n\to \infty}\frac{\#\{k\mid \lambda_{k+1}-\lambda_{k}=g, \, k\leq n\}}{n}=f.$$

\begin{theorem}[Bugeaud \cite{bug02}]\label{thm2}
Let $m \geq 2$ and let $q_m$ be the largest real root in $(1,2)$ of
\[
x^m - x^{m-1} - \cdots - x - 1.
\]
Let $\mathbf v_m = \mathbf v_{m,1}\mathbf v_{m,2}\cdots$ be the
$m$-bonacci word. 
Then, for all $n \geq 1$, one has
\[
\lambda_{n+1}(q_m)-\lambda_n(q_m)=d_m(j)
\quad \text{if and only if} \quad
\mathbf v_{m,n}=j,
\]
where
\[
d_m(1)=1,
\qquad
d_m(j)=q_m^{j-1}-q_m^{j-2}-\cdots - q_m - 1
\quad \text{for } j=2,\dots,m.
\]

Moreover, for $j=1,\dots,m$ the value $d_m(j)$ occurs in the gap sequence
$\bigl(\lambda_{n+1}(q_m)-\lambda_n(q_m)\bigr)_{n\geq 1}$
with frequency
$q_m^{-j}.$
\end{theorem}
For instance, in the golden mean case $m=2$, the gaps $1$ and $\varphi-1$ occur in $\{\lambda_{n+1}(\varphi)-\lambda_{n+1}(\varphi)\}$ with frequencies $1/\varphi$ and $1/\varphi^2$, respectively. 

\begin{example}[Fibonacci word]
The Fibonacci word corresponds to the case $m=2$. The first iterations of $\sigma_2$ on $1$ are
\begin{align*}
&\sigma_2(1)=12;\\
&\sigma_2^2(1)=\sigma_2(12)=\sigma_2(1)\sigma_2(2)=121;;\\
&\sigma_2^3(1)=12112;\\
&\sigma_2^4(1)=12112121\end{align*}
and for all $k$, $\sigma_2^k(1)$ is a prefix of $\mathbf v_2$. 
We have $d_2(1)=1$ and $d_2(2)=\varphi-1=1/\varphi$. 
The first three positive terms  of the gap sequence are 
$$\lambda_1(\varphi)-\lambda_0(\varphi)=d_2(\mathbf v_{2,1})=d_2(1)=1$$
and
$$\lambda_2(\varphi)-\lambda_1(\varphi)=d_2(\mathbf v_{2,2})=d_2(2)=\frac{1}{\varphi},\qquad \lambda_3(\varphi)-\lambda_2(\varphi)=d_2(\mathbf v_{2,3})=d_2(1)=1$$ Then 
$$\Lambda(\varphi)=\{0,1,1+1/\varphi,2+1/\varphi,\dots\}.$$


\end{example}

\begin{example}[Tribonacci word]
The first iterations of $\sigma_3$ on $1$ are
\begin{align*}
&\sigma_3(1)=12;\\
&\sigma_3^2(1)=\sigma_3(12)=\sigma_3(1)\sigma_3(2)=1213;\\
&\sigma_3^3(1)=1213121;\\
&\sigma_3^4(1)=1213121121312. \end{align*}
and for all $k$, $\sigma_3^k(1)$ is a prefix of $\mathbf v_3$. We have $\rho_3=\tau$, where $\tau$ is the greatest root of $x^3-x^2-x-1$, i.e., the Tribonacci constant. We have
$$d_3(1)=1\,\quad d_3(2)= \tau-1,   \quad d_3(3)=\tau^2-\tau-1.$$

The first three positive terms of the gap sequence are
\[
\lambda_1(\tau)-\lambda_0(\tau)
=
d_3(\mathbf v_{3,1})
=
d_3(1)
=
1,
\]
and
\[
\lambda_2(\tau)-\lambda_1(\tau)
=
d_3(\mathbf v_{3,2})
=
d_3(2)
=
\tau-1,
\qquad
\lambda_3(\tau)-\lambda_2(\tau)
=
d_3(\mathbf v_{3,3})
=
d_3(1)
=
1.
\]

Therefore
\[
\Lambda(\tau)
=
\{0,1,\tau,\tau+1,\dots\}.
\]
\end{example}

\subsection{Combinatorics on the $m$-bonacci word}
Theorem \ref{thm2} relates the occurrences of the gaps $d_m(j)$ in $\Lambda(q_m)$ to the digit distibution in the $m$-bonacci word. Therefore, it is useful in the sequel to state some definitions and results on the combinatorial properties of $\mathbf v_m$.

\begin{definition}[Weight and balanced word]
Let $w$ be a finite or infinite word over the alphabet $\mathcal A_m=\{1,\dots,m\}$.
For $j\in \mathcal A_m$, we denote by $|w|_j$ the number of occurrences of the letter $j$ in $w$, and call it the \emph{weight} of $j$ in $w$.

Let $b$ be a positive integer. An infinite word $\mathbf v$ over $\mathcal A_m$ is said to be $b$-\emph{balanced} for any two subwords $u$ and $w$ of $\mathbf v$ of the same length, and for every $j\in \mathcal A_m$, one has
\[
\bigl||u|_j - |w|_j\bigr| \leq b.
\]
\end{definition}
The next result describes the weights of the $m$-bonacci word. 
\begin{lemma}\label{l1bal}[\cite{balancetrib}]
For every prefix $v$ of the $m$-bonacci word $\mathbf v_m$ there exist $H$ and $x_0,x_1,\dots,x_H\in\{0,1\}$ such that 
\begin{equation}\label{prefix}(|v|_1,\dots,|v|_m)^T=\sum_{h=0}^H x_h M_m^h (1,0,\dots,0)^T \end{equation}
where $M_m\in \RR^m\times \RR^m$ is the incidence matrix of the Rauzy substitution $\sigma_m$, namely
\begin{equation}\label{Mm}M_{m}= 
\begin{pmatrix}
1 & 1 & 1 & \cdots & 1 & 1 \\
1 & 0 & 0 & \cdots & 0 & 0 \\
0 & 1 & 0 & \cdots & 0 & 0 \\
\vdots & \vdots & \vdots & \ddots & \vdots & \vdots \\
0 & 0 & 0 & \cdots & 1 & 0
\end{pmatrix}.
\end{equation}
Moreover, for every choiche of $H$ and $x_0,\cdots,x_H\in\{0,1\}$ there exists a prefix of $\mathbf v_m$ such that \eqref{prefix} holds. 
\end{lemma}

The $m$-bonacci words are balanced, the next result collects some estimates on the balance constants $b_m$ -- note that the case $m\geq 3$ makes use of the above result.

\begin{theorem}[\cite{bal,balancetrib,generalbalance}]\label{thm3}
For all $m\geq 2$ the $m$-bonacci word is $b$-balanced with 
\begin{enumerate}
\item $b=m-1$ for $m=2,3,4$ and this bound cannot be improved.
\item $b=\lceil \frac{m+1}{2}\rceil$ for $m=5,\dots,12$ and if $m=5$ this bound cannot be improved.
\end{enumerate}
Morevoer for all $m\geq 5$ the following estimate hold
$$b=\lceil \kappa m\rceil+12$$ 
where $\kappa:= \frac{2}{\pi} \int_{0}^{2\pi}
\frac{1-\cos x}{(5-4\cos x)\ln(5-4\cos x)} \, dx \approx 0.58.$
\end{theorem}

\subsection{The $m$-bonacci expansions}
Define the \emph{$m$-bonacci sequence} recursively
$$F^{(m)}_0=1, F^{(m)}_1=2, \cdots F^{(m)}_{m-1}=2^{m-1}$$
and for $n\geq m$
$$F^{(m)}_n=F^{(m)}_{n-1}+\cdots +F^{(m)}_{n-m}.$$
For $m=2$ one obtains the sequence $(F^{(2)})_n$ is $ 1,2,3,5,\dots$, the classical Fibonacci sequence (number A000045  in the OEIS classification) shifted by two indices. 
Similarly,  for $m=3$ the sequence $(F^{(3)}_n)= 1, 2, 4, 7,,\dots$ corresponds to the Tribonacci numbers (number A000073 in the OEIS classification), shifted by three indices.

In \cite{CSH72}, investigate the $m$-bonacci numeration system, also called the higher order Fibonacci expansions. Generalizing earlier results established in the Fibonacci case, the Zeckendorf representation \cite{Lek,zec}, a result in that paper states that every $N\in\NN$ admits an expansion of the form 
\begin{equation}\label{Nfib}N=\sum_{h=0}^{\infty}x_h F^{(m)}_h\qquad x_h\in\{0,1\}\end{equation}
where only finitely many digits $x_h$ are non zero and the digit sequence $(x_h)$ contains no blocks of $m$ consecutive $1$'s. Precisely, the greedy algorithm is the following.

Let $N\in\mathbb N$ and set $R_0:=N$.
For each $t\geq 0$, as long as $R_t>0$, define
\[
h_t:=\max\{h\geq 0:\ F_h^{(m)}\leq R_t\},
\]
and
\[
R_{t+1}:=R_t-F_{h_t}^{(m)}.
\]

Since $R_t$ is strictly decreasing, there exists $T\geq 0$ such that
$R_{T+1}=0$. The greedy $m$-bonacci expansion of $N$ is then
\[
N=\sum_{t=0}^{T}F_{h_t}^{(m)}
   =\sum_{h\geq 0}x_hF_h^{(m)},
\]
where
\[
x_h=
\begin{cases}
1,& \text{if } h=h_t \text{ for some } t,\\[2mm]
0,& \text{otherwise.}
\end{cases}
\]


Canonical representations are not unique in general. For instance, when $m=2$ we have 
$$3=F_0^{(2)}+F_1^{(2)}=F_2^{(2)}$$
so that $11(0)^\infty$ and $001(0)^\infty$ are both expansions of $3$.

\section{Gap conditions for the spectrum of $m$-bonacci numbers }\label{sec3}
In this section we prove our main results. We begin by proving Theorem \ref{thm1i} that is 

\begin{theorema}
Let $m\geq 2$ and $N\in \mathbb{N}$. Then there exists a constant
$\gamma_{m,N}>0$ such that
\[
\lambda_{n+N}(q_m)-\lambda_n(q_m)\geq N\,\gamma_{m,N}
\quad \text{for all } n\geq 1.
\]

The following choice of $\gamma_{m,N}$ satisfies the above inequality.
Let $(x_h)\in\{0,1\}^\infty$ be the canonical $m$-bonacci expansion of $N$ and let $b_m>0$ such that the $m$-bonacci word $\mathbf v_m$ is $b_m$-balanced.
Then
\begin{equation}\label{eqgamma}
\gamma_{m,N}
=\frac{1}{N}
\sum_{j=1}^m
\left(
\sum_{h\geq 0} x_h F^{(m)}_{h-j}
- b_m
\right)
d_m(j),
\end{equation}
where
 $d_m(j)=q_m^{j-1}-\sum_{k=0}^{j-2} q_m^k$ with $j=1,\dots,m$.
\end{theorema}
\begin{proof}
First of all we remark that, by Theorem \ref{thm2},
for all $m\geq 2$, $k,N \in \NN$
\begin{align*}\label{gapw}
    \lambda_{k+N}(q_m)-\lambda_{k}(q_m)&=
    \sum_{h=1}^N(\lambda_{k+h}(q_m)-\lambda_{k+h-1}(q_m))\\
    &=\sum_{j=1}^m |\{h\mid \lambda_{k+h}(q_m)-\lambda_{k+h-1}=d_m(j),\,h=1,\dots,N\}|d_m(j)\\
    &=\sum_{j=1}^m |\mathbf v_{k+1}\cdots \mathbf v_{k+N}|_j d_m(j).
\end{align*}
As $\mathbf v_m$ is a balanced word it admits a lowest balance constant $b_m$, see Theorem \ref{thm3}. We then deduce from the above equality the lower estimate 
\begin{equation}\label{gapw}
    \lambda_{k+N}(q_m)-\lambda_{k}(q_m) \geq \sum_{j=1}^m\left(|\mathbf v_{m,1}\cdots \mathbf v_{m,N}|_j-b_m\right)d_m(j)\end{equation}
for all $m\geq 2$, $k,N \in \NN$.

Let $(x_h)$ be the canonical $m$-bonacci expansion of $N$.
By the converse statement of Lemma \ref{l1bal}, there exists a
prefix $v$ of $\mathbf v_m$ such that

\[
(|v|_1,\dots,|v|_m)^T
=
\sum_{h\ge0}x_h M_m^h(1,0,\dots,0)^T .
\]

Now, a simple inductive argument (already remarked in \cite{generalbalance} with a slightly different notation) shows that if $M_m$ is the incidence matrix of $\sigma_m$, defined in \eqref{Mm} then
$$M^h(1,0,\dots,0)^T=(F^{(m)}_{h-1},F^{(m)}_{h-2},\dots, F^{(m)}_{h-m})^T$$
where $F^{(m)}_{-1}=1$ and $F^{(m)}_{-k}=0$ for all $k> 1$.
Then we obtain
\[
|v|_j
=
\sum_{h\ge0}x_hF_{h-j}^{(m)}
\qquad (j=1,\dots,m).
\]

Moreover,

\[
|v|
=
\sum_{j=1}^m|v|_j
=
\sum_{h\ge0}x_h
\sum_{j=1}^mF_{h-j}^{(m)}
=
\sum_{h\ge0}x_hF_h^{(m)}
=
N.
\]

Therefore $v$ is the prefix
$\mathbf v_{m,1}\cdots\mathbf v_{m,N}$ and hence

\[
|\mathbf v_{m,1}\cdots\mathbf v_{m,N}|_j
=
\sum_{h\ge0}x_hF_{h-j}^{(m)}.
\]
This, together with \eqref{gapw}, concludes the proof of the theorem.


\end{proof}

The next results discuss the cases $m=2,3$.

\begin{corollary}[Fibonacci case]\label{cor1}
For all $N\in \NN$ we have the gap condition
\begin{equation}\label{gengapfib}
\lambda_{2,N+k}-\lambda_{2,k}\geq N\gamma_{2,N} \quad \text{for all } k\in \ZZ,
\end{equation}
where
$$\gamma_{2,N}:=\frac{1}{\varphi}+\frac{1}{\varphi^3} -\frac{1}{N}\left(\frac{1}{\varphi^4}+\varphi\right).$$
The estimate is sharpened in the case $N=2$ by setting  $$\gamma_{2,2}=\frac{\varphi}{2}.$$

\end{corollary}
\begin{proof}
Set for brevity $F^{(2)}_n=F_n$ (so that $F_0=1, F_1=2,F_2=3\cdots$, in particular $F_n$ is the $n+2$th Fibonacci number) and let $(x_h)$ be the canonical Fibonacci expansion of $N$. Note that $d_2(1)=1$ and $d_2(2)=1/\varphi$. Since the balance constant $b_2$ is equal to $1$, by Theorem \ref{thm1i} we have that 
\eqref{gengapfib} holds with 
\begin{align*}\gamma_{2,N}
&=
\frac{1}{N}
\sum_{j=1}^2
\left(
\sum_{h=0}^\infty x_h F_{h-j} - 1
\right)d_2(j)\\
&=
\frac{1}{N}
\left[
\left(\sum_{h=0}^\infty x_h F_{h-1} -1\right)
+
\left(\sum_{h=0}^\infty x_h F_{h-2} -1\right)\frac{1}{\varphi}
\right] \\
&=
\frac{1}{N}
\left(
\sum_{h=0}^\infty x_h\left(F_{h-1}+F_{h-2}\frac{1}{\varphi}\right)
-\varphi
\right).
\end{align*}
Since 
$$F_{h-1}=\frac{1}{\varphi} F_h +  \left(-\frac{1}{\varphi}\right)^{h+2};\quad F_{h-2}=\frac{1}{\varphi^2} F_h -  \left(-\frac{1}{\varphi}\right)^{h+2} $$
we obtain
\begin{align*}\gamma_{2,N}&=\frac{1}{ N}\left(\left(\frac{1}{\varphi}+\frac{1}{\varphi^3}\right) \sum_{h=0}^\infty x_h F_h+ \frac{1}{\varphi^4}\sum_{h=0}^\infty  \frac{x_h}{(-\varphi)^h}-\varphi\right)\\
&=\frac{1}{\varphi}+\frac{1}{\varphi^3} +\frac{1}{ N}\left(\frac{1}{\varphi^4}\sum_{h=0}^\infty \frac{x_h}{(-\varphi)^h}-\varphi\right)\\
&\geq\frac{1}{\varphi}+\frac{1}{\varphi^3}  -\frac{1}{ N}\left(\frac{1}{\varphi^3}\sum_{h=1}^\infty \frac{1}{\varphi^{2h}}+\varphi\right)\\
&= \frac{1}{\varphi}+\frac{1}{\varphi^3} -\frac{1}{N}\left(\frac{1}{\varphi^4}+\varphi\right).
\end{align*}

In the particular case $N=2$, we note that the subword $22$ does not occur
in the Fibonacci word.
Therefore,
\[
\lambda_{2,k+2}-\lambda_{2,k}\ge d_2(1)+d_2(2)=\varphi,
\]
and one can choose $\gamma_{2,2}=\varphi/2$.

\end{proof}
\begin{corollary}[ Tribonacci case]\label{cor2}
Let $N\in \NN$ and let $(x_h)$ be the canonical Tribonacci expansion of $N$. Then we have the gap condition
\begin{equation}\label{gengaptr}
\lambda_{3,N+k}-\lambda_{3,k}\geq N\gamma_{3,N} \quad \text{for all } k\in \ZZ,
\end{equation}
where
\begin{align*}\gamma_{3,N}&:=\frac{1}{N}\left(\sum_{h=0}^\infty  x_h\tau^{h}-2(\tau^2-1)\right)\\&\geq \frac{1}{C_1\tau^3} -\frac{1}{N}\left({\frac{1}{2C_1\tau^3}\log_\tau \left( \frac{N +\frac{1}{2}}{C_1} \right)+2\tau^2-2-\frac{1}{C_1\tau^3})}\right).\end{align*}
where $\tau$ is the Tribonacci constant (i.e. the greatest root of $\tau^3-\tau^2-\tau-1$) and  $$C_1:= \frac{\tau^3}{4\tau^2 - 3\tau - 3}\sim 0.18. $$

The estimate is sharpened in the case $N=2$ by setting  $$\gamma_{3,2}=\frac{\tau^2-\tau}{2}.$$
\end{corollary}

\begin{proof}
Set for brevity $F^{(3)}_n=T_n$ (so that $T_0=1, T_1=2,T_2=4,T_3=7 \cdots$, in particular $T_n$ is the $n+3$th Tribonacci number) and let $(x_h)$ be the canonical $m$-bonacci (i.e., Tribonacci) expansion of $N$. Note that $b_3=2$ and

$$d_2(1)=1\,\quad d_2(2)= \tau-1=\frac{1}{\tau}+\frac{1}{\tau^2},   \quad d_2(3)=\tau^2-\tau-1=\frac{1}{\tau}.$$

By Theorem \ref{thm1i} we have that \eqref{gengaptr} holds with 
\begin{equation}\label{r}\begin{split}\gamma_{3,N}=&\frac{1}{N}\sum_{h=0}^\infty x_h\left(T_{h-1}+\frac{\tau+1}{\tau^2}T_{h-2}+\frac{1}{\tau}T_{h-3}\right)-2\frac{\tau^2-1}{N}.\end{split}\end{equation}
One can check by induction  that 
\begin{equation}\label{rec} \tau^h= T_{h-1}+\frac{\tau+1}{\tau^2} T_{h-2}+ \frac{1}{\tau}T_{h-3}.\end{equation}
Replacing the above expression in \eqref{r}  we have that 
\begin{align}\gamma_{3,N}=\frac{1}{N}\left( \sum_{h=0}^\infty x_h \tau^h-2(\tau^2-1)\right).\label{g}
\end{align}
By \cite{binet} and by the fact  that $T_h$ is the $h+3$-th Tribonacci number, we have  \begin{equation}\label{binet} |T_{h}-C_1\tau^{h+3}|\leq \frac{1}{2} \end{equation}
for all $h$ and, consequently, the lower estimate
$$\tau^h\geq \frac{1}{C_1\tau^3}\left(T_h-\frac{1}{2}\right)$$
yielding
\begin{align}\gamma_{3,N}\geq\frac{1}{C_1\tau^3}-\frac{1}{N}\left( \frac{1}{2C_1\tau^3}\sum_{h=0}^\infty x_h +2(\tau^2-1)\right).\label{g1}
\end{align}
To estimate the sum $\sum_{h=0}^\infty x_h$, let $H=H(N)$ be such that $N\in[T_H,T_{H+1})$ so that $x_h=0$ for all $h>H$.  
Then $N\geq T_{H}$ and by \eqref{binet} we deduce
$$N\geq C_1\tau^{H+3} - \frac{1}{2}$$
and
$$\sum_{h=0}^\infty x_h\leq H+1 \le \log_\tau\left(\frac{N+1/2}{C_1}\right)-2$$
Plugging the above inequality in \eqref{g1} we conclude
\begin{align*}\gamma_{3,N}&\geq
\frac{1}{C_1\tau^3} -\frac{1}{N}\left( \frac{1}{2C_1\tau^3}\log_\tau \left( \frac{N +\frac{1}{2}}{C_1} \right)+2\tau^2-2-\frac{1}{C_1\tau^3})\right).
\end{align*}

In the particular case $N=2$, we note that the subwords of length 2 in $\mathbf v_3$ are the following $\{11,12,13,21,31\}$, see \cite{seminalrauzy}. We then have
$$\lambda_{k+2}(\tau)-\lambda_k(\tau)\geq d_3(1)+d_3(3)=\tau^2-\tau.$$
\end{proof}

\section{Some remarks on the density of $\Lambda(q_m)$}\label{sec4}

The gap frequency of the spectrum $\Lambda(q)$ is deeply related to its \emph{upper density} \cite{beu}

\begin{equation*}
D^+(\Lambda(q)) := \lim_{R\to +\infty} \frac{n^+(R)}{R}
\end{equation*}
where $n^+(R)$ denotes the largest number of elements of $\Lambda(q)$ contained in an interval of length $R$, see \cite{komlorbook} for the existence of the limit. In \cite{FW02} Feng and Weng proved that every Pisot spectrum is associated with a substitution. When the substitution has the further property of being \emph{primitive} -- namely the power $M^k$ of its incidence matrix $M$ has positive entries for all sufficiently large $k$--  the associated symbolic dyanamics is strictly ergodic, from which follows that the gaps admit \emph{uniform frequencies}. Therefore the average gap length is independent of the position of the observation window and the upper density is given by
\begin{equation}\label{density}
D^+(\Lambda(q)) =\frac{1}{\sum_{k=1}^{N(q)} f_k(q)g_k(q)}.
\end{equation}
where $(f_k(q),g_k(q))$ are the finite, say $N(q)$, (uniform) frequency-gap couples associated to $\Lambda(q)$. In general, it is not known whether every Pisot number yield a uniform frequency of the gaps in the its spectrum, but several classes of Pisot numbers, including the $m$-bonacci numbers, with this property were provided.  
For instance, in the golden mean case $m=2$, such couples are $(\frac{1}{\varphi},1), (\frac{1}{\varphi^2},\varphi-1)$, $N(q_2)=2$ and the upper density is $D^+(\Lambda(\varphi))=\frac{1}{\frac{1}{\varphi}+\frac{1}{\varphi^3}}=\frac{\varphi^2}{\sqrt{5}}$ -- this can also be deduced by the fact that every sufficienly large $\lambda$ belongs to $\Lambda(\varphi)$ if and only if $ \lambda/\varphi$ belongs to the celebrated Fibonacci quasicrystal \cite{H04,MPP15}, whose density is $\frac{\varphi}{\sqrt{5}}$. 

 In \cite{GH06} the following result is proved. Let $q$ be the Pisot root of 
$$x^m-ax^{m-1}-\cdots ax -b$$
where $a\geq b\geq 1$. Let $r=\lfloor q\rfloor=a$. Then the gaps in $\Lambda^{r}(q)$ are 
$$ 1,\quad q-a, \quad q^2-a q-a,\quad\dots \quad  q^{m-1}-a q^{m-2}-\cdots a q-a$$ and they respectively 
occur with uniform frequency 
$$\frac{q^{m-1}}{Q},\,\frac{q^{m-2}}{Q},\dots, \frac{1}{Q}$$
where $Q:=q^{m-1}+q^{m-2}+\cdots+q+1=\frac{q^m-1}{q-1}$. This allows to compute the upper density 
\begin{align*}D^+(\Lambda^r(q))&=
\frac{Q}
{mq^{m-1}-a\sum_{k=0}^{m-2}(k+1)q^{k}}\\
&=\frac{Q}{mq^{m-1}-a\frac{(m-1)q^{m}-mq^{m-1}+1}{(q-1)^2}}\\
&=\frac{(q^m-1)(q-1)}{mq^{m+1}-(2m+am-a)q^{m}+m(a+1)q^{m-1}-a}.
\end{align*}
In the case of $m$-bonacci Pisot numbers, the above formula yields 
$$D^+(\Lambda(q_m))= \frac{1}{\frac{1}{q_m}+\sum_{k=2}^m \frac{q_m^{k-1}-q_m^{k-2}-\cdots-q_m-1}{q_m^k}}=\frac{(q^m_m-1)(q_m-1)}{mq_m^{m+1}-(3m-1)q_m^{m}+2mq_m^{m-1}-1}.$$

The density provides an average distribution of the gaps. 
\medskip
 By \cite{bkl98}, see also \cite{komlorbook}, since $\Lambda(q_m)$ has finite upper density, then there exists $\gamma_{N,m}$ satisfying  \eqref{gamma} and 
\begin{equation}\label{denlim}\lim_{N\to \infty}\gamma_{N,m}=\frac{1}{D^+(\Lambda(q_m))}=\frac{1}{q_m}+\sum_{k=2}^m \frac{q_m^{k-1}-q_m^{k-2}-\cdots-q_m-1}{q_m^k}.\end{equation}
A trivial choice for $\gamma_{N,m}$ is the the smallest gap in $\Lambda(q_m)$, that is $\frac{1}{q_m^{m-1}}$, because 
$$\lambda_{k+N}(q_m)-\lambda_{k}(q_m)=\sum_{h=0}^{N-1}\lambda_{k+h+1}(q_m)-\lambda_{k+h}(q_m)\geq \frac{N}{q_m^{m-1}}.$$
Such an estimate is far from the “ideal” choice of $\gamma_{N,m}$ suggested by \eqref{denlim}, since it only reflects the minimal gap and ignores the finer distribution of gaps in the spectrum. An explicit estimate consistent with \eqref{denlim} for $m=2$ is for instance the estimate provided in \eqref{lambdaest}. 

\section{Conclusions}\label{sec5}
In this paper we established  explicit gap estimates for the spectrum of the  $m$-bonacci Pisot numbers. Our methodology allows for a finer investigation of the structure of the spectra and it relies on the combinatorical properties of the underlying $m$-bonacci words, particularly the existence of forbidden subwords, the balance property (factors of the same lengths have bounded variations in their digit distributions), and the digital representations of integers via $m$-bonacci number systems.  We then applied this general framework to derive further estimates for the Fibonacci and Tribonacci cases. Some considerations on the density of the spectra conclude our investigation.

\bibliographystyle{abbrv}
\bibliography{rauzy}
\end{document}